\documentclass[12pt,reqno]{amsart}
\usepackage{amsmath, amsthm, amssymb}
\makeatletter
\@namedef{subjclassname@2010}{%
  \textup{2010} Mathematics Subject Classification}
\makeatother

\usepackage{color}
\usepackage{comment}
\usepackage{amscd}

\advance \topmargin by -\headheight
\advance \topmargin by -\headsep

\evensidemargin \oddsidemargin
\DeclareFontFamily{U}{wncy}{}
\DeclareFontShape{U}{wncy}{m}{n}{%
   <5>wncyr5%
   <6>wncyr6%
   <7>wncyr7%
   <8>wncyr8%
   <9>wncyr9%
   <10>wncyr10%
   <11>wncyr10%
   <12>wncyr6%
   <14>wncyr7%
   <17>wncyr8%
   <20>wncyr10%
   <25>wncyr10}{}
\DeclareMathAlphabet{\cyr}{U}{wncy}{m}{n}

\newtheorem{theorem}{Theorem}[section]
\newtheorem{lemma}[theorem]{Lemma}
\newtheorem{cor}[theorem]{Corollary}

\newtheorem{prop}[theorem]{Proposition}

\theoremstyle{definition}

\theoremstyle{remark}

\numberwithin{equation}{section}

\def\bft{{\mathbf t}}

\def\bfx{{\mathbf x}}
\def\bfy{{\mathbf y}}
\def\bfz{{\mathbf z}}
\def\AA{{\mathbf A}}

\def\calU{{\mathcal U}}

\def\calX{{\mathcal X}}

\def\A{{\mathbb A}}

\def\G{{\mathbb G}}
\def\F{{\mathbb F}}\def\P{{\mathbb P}}
\def\R{{\mathbb R}}
\def\Z{{\mathbb Z}}\def\Q{{\mathbb Q}}

\def\ov{\overline}
\def\lra{\longrightarrow}

\specialcomment{comblue}{\begin{comblue}}{\end{comblue}}
\excludecomment{comblue}

\specialcomment{comgreen}{\begin{comgreen}}{\end{comgreen}}
\excludecomment{comgreen}

\def\deg{{\rm deg}}

\def\Br{{\rm Br}}

\def\Gal{{\rm Gal}}
\def\Pic{{\rm Pic}}

\def\val{{\rm val}}
\def\Brv{{{\rm Br}_{\rm vert}}}
\def\inv{{\rm inv}}

\begin{document}

\title{The Hardy--Littlewood conjecture and
rational points}

\date{\today}

\author{Yonatan Harpaz}
\address{Institute for Mathematics, Astrophysics and Particle Physics\\
Radboud University Nijmegen\\Heyendaalseweg 135\\
6525 AJ Nijmegen\\the Netherlands }
\email{harpazy@gmail.com}

\author{Alexei N. Skorobogatov}
\address{Department of Mathematics\\ Imperial College London \\ 
SW7 2BZ\\U.K.}
\address{Institute for the Information Transmission Problems\\ 
Russian Academy of Sciences\\
19 Bolshoi Karetnyi\\ Moscow\\ 127994\\ Russia}
\email{a.skorobogatov@imperial.ac.uk}

\author{Olivier Wittenberg}
\address{D\'epartement de math\'ematiques et applications\\ \'Ecole normale sup\'erieure\\ 45 rue
d'Ulm\\ 75230\\ Paris Cedex 05\\ France}
\email{wittenberg@dma.ens.fr}



\maketitle

\begin{abstract}
Schinzel's Hypothesis (H) was used by Colliot-Th\'el\`ene and Sansuc, 
and later by Serre, Swinnerton-Dyer and others, to prove that 
the Brauer--Manin obstruction controls the Hasse principle and 
weak approximation on pencils of conics and similar varieties. 
We show that when the ground field is $\Q$ and the degenerate geometric 
fibres of the pencil are all defined over $\Q$, one can use this method 
to obtain unconditional results by replacing Hypothesis (H) with 
the finite complexity case of the generalised Hardy--Littlewood 
conjecture recently established by Green, Tao and Ziegler. 
\end{abstract}

\section*{Introduction}

The finite complexity case of the generalised Hardy--Littlewood conjecture
recently proved by Green and Tao \cite{GT1, GT2} and Green--Tao--Ziegler
\cite{GTZ} is of fundamental importance to number theory.
The aim of this note is to explore some of its consequences
for the Hasse principle and weak approximation
on algebraic varieties over $\Q$.

Hasse used Dirichlet's theorem on primes in an arithmetic
progression to deduce what is now called the Hasse principle
for quadratic forms in four variables from the global
reciprocity law and the Hasse principle
for quadratic forms in three variables, itself a corollary
of global class field theory (see \cite{Hasse}, p.~16 and p.~87).
This `fibration method' was taken up
by Colliot-Th\'el\`ene and Sansuc in \cite{CSan}. They showed
that Schinzel's Hypothesis~(H), a vast generalisation of
Dirichlet's theorem in which $at+b$ is replaced by a finite collection of
irreducible polynomials in $t$ of arbitrary degree \cite{SchSi58}, 
has strong consequences
for the Hasse principle and weak approximation on conic bundles
and some other pencils of varieties over $\Q$.
The method of \cite{CSan} was extended to number fields and 
generalised in various directions by Serre (see \cite{Serre}, Ch. II, Annexe), 
Swinnerton-Dyer (see \cite{Sw94}, \cite{SwD} and references in that paper), 
Colliot-Th\'el\`ene and others \cite{CS,CSS,W,Wei}.
In this note we show that when the ground field is $\Q$ and
the degenerate geometric fibres of the pencil are all
defined over $\Q$, one can apply the methods of 
the aforementioned papers to obtain unconditional results by using
the theorem of Green, Tao and Ziegler in place of Schinzel's Hypothesis~(H).

Additive combinatorics was recently applied to the study of rational points
by Browning, Matthiesen and one of the authors in \cite{BMS}. 
The approach of \cite{BMS} uses the descent method of Colliot-Th\'el\`ene 
and Sansuc; it crucially relies on
the work of Matthiesen on representation of linear polynomials by
binary quadratic forms \cite{LM,LM1} in order to prove the Hasse principle
and weak approximation for certain varieties appearing after descent.
In this paper we give a short proof of most of the results of \cite{BMS}
as well as a generalisation to certain equations involving norms of cyclic (and some
non-cyclic) extensions. Our approach is
based directly on the Green--Tao--Ziegler theorem and on various generalisations
of the method of Hasse, Colliot-Th\'el\`ene--Sansuc and Swinnerton-Dyer
mentioned above, thus avoiding the use of descent and the results of \cite{LM,LM1}.

We recall the theorem of Green, Tao and Ziegler and give its first
corollaries in Section \ref{one}. In Section \ref{Hyp} we compare
Schinzel's Hypothesis~(H), in the form of Hypothesis~${\rm (H_1\mkern-1.5mu)}$,
with Proposition \ref{p2}, a consequence of the Green--Tao--Ziegler theorem.
We prove our main results and deduce their first applications in 
Section \ref{fib}. 
In Section \ref{four} we consider applications to representation of norms
(and products of norms) by products of linear polynomials.

This paper was started during the special semester `Rational points
and algebraic cycles' at the {\em Centre Interfacultaire Bernoulli} 
of the \'Ecole Polytechnique F\'ed\'erale de Lausanne, which we would like
to thank for its hospitality. The second author would also like to thank
the Hausdorff Research Institute for Mathematics in Bonn for its
hospitality during the `Arithmetic and geometry' program.
The authors are grateful to T.D. Browning, J.-L. Colliot-Th\'el\`ene 
and D. Schindler for their comments, and to
B. Viray for pointing out a mistake in the first version
of this paper.

\section{A corollary of the generalised Hardy--Littlewood conjecture} \label{one}

In a series of papers Green and Tao \cite{GT1, GT2} and
Green--Tao--Ziegler \cite{GTZ} proved the generalised Hardy--Littlewood
conjecture in the finite complexity case. The following
qualitative statement is \cite[Cor.~1.9]{GT1}.

\begin{theorem}[Green, Tao, Ziegler] \label{gtz}
Let $L_1(x,y),\ldots,L_r(x,y)\in\Z[x,y]$ be pairwise non-proportional
linear forms, and let $c_1,\ldots,c_r\in \Z$.
Assume that for each prime $p$, there exists $(m,n)\in\Z^2$
such that~$p$ does not divide $L_i(m,n)+c_i$ for any $i=1,\ldots,r$.
Let $K\subset\R^2$ be an open convex cone containing a point
$(m,n)\in\Z^2$
such that $L_i(m,n)>0$ for $i=1,\ldots,r$.
Then there exist infinitely
many pairs $(m,n)\in K\cap\Z^2$ such that $L_i(m,n)+c_i$
are all prime.
\end{theorem}

\medskip

We shall use the following easy corollary of this result.
For a finite set of rational primes $S$ we write $\Z_S=\Z[S^{-1}]$.

\begin{prop} \label{p1}
Suppose that we are given $(\lambda_p,\mu_p)\in\Q_p^2$
for $p$ in a finite set of primes $S$, and a positive real constant $C$.
Let $e_1,\ldots, e_r$ be 
elements of $\Z_S$. Then there exist infinitely many pairs
$(\lambda,\mu)\in\Z_S^2$ such that
\begin{enumerate}
\item
$\lambda>C\mu>0$;

\item
$(\lambda,\mu)$ is close to $(\lambda_p,\mu_p)$ in the $p$-adic topology
for $p\in S$;
\item
$\lambda-e_i\mu=p_i u_i$ with $u_i\in \Z_S^*$, for $i=1,\ldots,r$,
where $p_1,\ldots, p_r$ are primes not in~$S$ 
such that $p_i=p_j$ if and only if $e_i=e_j$.
\end{enumerate}
\end{prop}
{\em Proof.} By eliminating repetitions we can assume that
$e_1,\ldots, e_r$ are pairwise different.
We can multiply $\lambda,\mu$ and all $\lambda_p,\mu_p$ by a product of
powers of primes from $S$, and so assume without loss of generality
that $(\lambda_p,\mu_p)\in\Z_p^2$ for $p\in S$.
Replacing $C$ with a larger constant we assume $C>e_i$, $i=1,\ldots,r$.
Using the Chinese remainder theorem, we find
$\lambda_0 \in \Z$ such that $\lambda_0-\lambda_p$ is divisible by
a sufficiently high power $p^{n_p}$ for all $p\in S$, and similarly for
$\mu_0 \in \Z$. In doing so we can assume that $\lambda_0 > C \mu_0 > 0$,
in particular, $\lambda_0-e_i\mu_0>0$ for all $i$.

Let $d$ be a product of powers of primes from $S$ such that
$de_i\in\Z$ for \mbox{$i=1,\ldots, r$}.
Let us write $d(\lambda_0-e_i\mu_0)=M_i c_i$,
where $M_i$ is a product of powers
of primes from $S$, and $c_i\in\Z$ is coprime to the primes in $S$.
Let $N$ be a product of primes in $S$ such that $N>c_i-c_{j}$ for
any $i$ and $j$.
Let $M = \prod_{p\in S}p^{m_p}$ where
$$ m_p \geq \max\{n_p,\val_p(N)+\val_p(M_i)\}, \quad i=1,\ldots,r. $$
Then $N$ divides $M/M_i$ for each $i$.
We now look for $\lambda$ and $\mu$ of the form
\begin{equation} \label{E3}
\lambda = \lambda_0 + Mm, \quad \mu = \mu_0 + Mn, \quad (m,n) \in \Z^2.
\end{equation}
Write $L_i(x,y)=M_i^{-1}Md(x-e_iy)$, then
\begin{equation} \label{EE3}
\lambda - e_i\mu = d^{-1}M_i(L_i(m,n)+c_i)
\end{equation}
for each $i=1,\ldots, r$. Let us check that the linear functions
$L_i(x,y) + c_i$ satisfy the
condition of Theorem~\ref{gtz}. For $p\in S$, the integer
$L_i(0,0) + c_i$ is non-zero modulo $p$ for each $i$.
Now let $p$ be a prime not in $S$.
Since the determinant of the homogeneous part of
the affine transformation (\ref{E3}) is in $\Z_S^*$ and each $M_i^{-1}d(\lambda-e_i\mu)$
equals $M_i^{-1}d \in \Z_S^*$ at the point $(\lambda,\mu)=(1,0)$, 
we see that there is $(m,n)\in\Z^2$
such that $\prod_{i=1}^r(L_i(m,n) + c_i)$ is not divisible by $p$.

We now choose an open convex cone $K$. Choose $(m_0,n_0)\in\Z^2$,
$m_0>Cn_0>0$, for which the positive integers $L_i(m_0,n_0)$ are pairwise different. 
After re-ordering the subscripts, we see that the inequalities
$$x>Cy>0, \quad L_1(x,y)>\ldots >L_r(x,y)>0$$
hold for $(x,y)=(m_0,n_0)$.
Define $K\subset\R^2$ by these inequalities.
We can apply Theorem~\ref{gtz} to the linear functions
$L_i(x,y)+c_i$ and the cone $K$. Thus there exist infinitely many
pairs $(m,n)\in K\cap \Z^2$ such that $L_i(m,n)+c_i=p_i$,
where $p_i$ is a prime not in $S$, for $i=1,\ldots,r$.
The coefficients of each $L_i(x,y)$ are divisible by $N$, hence
$$L_i(m,n)-L_{i+1}(m,n)\geq N> c_{i+1}-c_i.$$
Thus $p_i>p_{i+1}$ for each $i=1,\ldots,r-1$, so all the primes $p_i$
are pairwise different.
Since $n>0$ and $m>Cn$ we see that $\mu=\mu_0+Mn>0$
and $\lambda=\lambda_0+Mm>C\mu$.

By (\ref{EE3}) we see that
$\lambda-e_i\mu$ differs from $L_i(x,y)+c_i$
by an element of $\Z_S^*$, so the proof is now complete. $\Box$

\medskip

Proposition~\ref{p1} can be used to study the Hasse principle and 
weak approximation for rational points. In
the proof of the following result, which is modelled on the original proof 
of Hasse (in the version of \cite[Prop.~2]{CSan}), it replaces
Dirichlet's theorem on primes in an arithmetic progression. 

For a field extension $K/\Q$ of degree $n$ we denote by 
$N_{K/\Q}(\bfx)$ the corresponding norm form of degree $n$
in $n$ variables $\bfx=(x_1,\ldots,x_n)$,
defined by choosing a basis of the $\Q$-vector space $K$.

\begin{theorem} \label{t1}
Let $K_i$ be a cyclic extension of $\Q$ of degree $d_i$ and
let $b_i\in\Q^*$, $e_i\in\Q$, for $i=1,\ldots,r$.
Then the affine variety $V\subset\A^2\times \A^{d_1}\times\ldots\times\A^{d_r}$
over $\Q$ defined by
\begin{equation}
b_i(u-e_iv)=N_{K_i/\Q}(\bfx_i)\not=0, \quad i=1,\ldots,r, \label{E1}
\end{equation}
satisfies the Hasse principle and weak approximation.
\end{theorem}
{\em Proof.} We are given $M_p\in V(\Q_p)$ for each prime $p$, and
$M_0\in V(\R)$. Let $S$ be the set of places of $\Q$ where we need
to approximate. We include the real place in $S$. 
Note that the set of real points $(u,v,\bfx_1,...,\bfx_r) \in V(\R)$ 
for which $(u,v)\in \Q^2$ is dense in $V(\R)$, and so it will be 
enough to prove the claim in the case when the coordinates $u$ and $v$
of $M_0$ are in $\Q$. By a $\Q$-linear change of variables 
we can assume without loss of generality that $M_0$ 
has coordinates $(u,v)=(1,0)$. Then we have $b_i>0$ 
whenever $K_i$ is totally imaginary.

We enlarge $S$ so that $b_i\in\Z_S^*$, $e_i\in\Z_S$ and the field
$K_i$ is not ramified outside~$S$, for all $i=1,\ldots,r$. Thus
for each $p\in S$ we now have a pair $(\lambda_p,\mu_p)\in\Q_p^2$
such that
$$b_i(\lambda_p-e_i\mu_p)=N_{K_i/\Q}(\bfx_{i,p})\not=0,\quad i=1,\ldots,r,$$
for some $\bfx_{i,p}\in (\Q_p)^{d_i}$.
Let $C$ be a large positive constant to be specified later, such that
$C>e_i$ for $i=1,\ldots, r$.
An application of Proposition~\ref{p1} produces 
$(\lambda,\mu)\in\Z_S^2$, $\lambda>C\mu>0$,
such that for each $i$ the number $b_i(\lambda-e_i\mu)$
is a local norm for $K_i/\Q$ at each finite place of $S$.
This is also true for the real place because
$b_i>0$ whenever $K_i$ is totally imaginary, and
$\lambda-e_i\mu>0$ for all $i$. Moreover, for each $i$
we have $b_i(\lambda-e_i\mu)=p_iu_i$, where $p_i$ is a prime not in $S$
and $u_i\in\Z_S^*$. Recall that $p_i=p_j$ if and only if $e_i=e_j$.

Let $(K_i/\Q,b_i(\lambda-e_i\mu))\in\Br(\Q)$ be the class
of the corresponding cyclic algebra. By continuity we have
$\inv_{p}(K_i/\Q,b_i(\lambda-e_i\mu))=0$
for any $p\in S$, and also
$\inv_{\R}(K_i/\Q,b_i(\lambda-e_i\mu))=0$.
Next, $b_i(\lambda-e_i\mu)$ is a unit at every prime $p\not\in S\cup\{p_i\}$,
hence we obtain
$$\inv_{p}(K_i/\Q,b_i(\lambda-e_i\mu))=0$$
for any $p\not=p_i$. The global reciprocity law now implies
$$\inv_{p_i}(K_i/\Q,b_i(\lambda-e_i\mu))=\inv_{p_i}(K_i/\Q,p_i)=0,$$
and since $K_i/\Q$ is unramified outside $S$,
the prime $p_i$ splits completely in $K_i$. In particular, $b_i(\lambda-e_i\mu)$
is a local norm at every place of $\Q$. By Hasse's theorem it is
a global norm, so that $$b_i(\lambda-e_i\mu)=N_{K_i/\Q}(\bfx_i)\not=0$$
for some $\bfx_i\in \Q^{d_i}$. This proves the Hasse principle for $V$.

Let us now prove weak approximation. Write $d=d_1\ldots d_r$.
Using weak approximation in $\Q$ we find a positive rational
number $\rho$ that is
$p$-adically close to $1$ for each prime $p\in S$,
and $\rho^d$ is close to $\lambda>0$ in the real topology. We now make
the change of variables
$$\lambda=\rho^d\lambda', \quad \mu=\rho^d\mu',\quad
\bfx_i=\rho^{d/d_i}\bfx'_i,\ \ i=1,\ldots, r.$$ Then
$(\lambda',\mu')$ is still close
to $(\lambda_p,\mu_p)$ in the $p$-adic topology for each prime $p\in S$.
In the real topology $(\lambda',\mu')$ is close to $(1,\mu/\lambda)$.
Since $0<\mu/\lambda<C^{-1}$, by choosing a large enough $C$
we ensure that $(\lambda',\mu')$ is
close to $(1,0)$. We can conclude by using weak approximation
in the norm tori $N_{K_i/\Q}(\bfz)=1$. $\Box$

\medskip

\noindent{\bf Remarks}. (1) In the case when all the fields $K_i$ are quadratic
this result is Thm.~1.2 of \cite{BMS}.
Eliminating $u$ and $v$ one sees that $V$ is then isomorphic to
an open subset of a complete intersection of $r-2$ quadrics in $\A^{2r}$
of a very special kind.
These intersections of quadrics are important because of their relation to conic bundles, first
pointed out by Salberger in \cite{Sa86}. In \cite{CSan2} Colliot-Th\'el\`ene and Sansuc 
proved the main results of descent theory and gave
a description of universal torsors by explicit equations.
Given the conclusion of Theorem \ref{t1} their results imply that
the Brauer--Manin obstruction is the only
obstruction to weak approximation on conic
bundles over $\P^1_\Q$ such that all the degenerate fibres are over $\Q$-points
(this is Thm.~1.1 of \cite{BMS};
see \cite{CSan2}, Thm.~2.6.4 and Chapter III, or \cite{Sk}, Prop.~4.4.8 and Cor.~6.1.3).

(2) In Section \ref{nn} below we give a somewhat different proof of Theorem \ref{t1}
deducing it from the main result of this paper, Theorem~\ref{t2}.

\section{Comparison with Schinzel's Hypothesis~(H)} \label{Hyp}

Our aim in this section is to show that in the study of 
rational points on pencils of varieties Proposition~\ref{p1} can replace
Schinzel's Hypothesis~(H). In such applications it is often more 
convenient to use a consequence of~(H),
stated in \cite[Section 4]{CS} under the name of 
Hypothesis~${\rm (H_1\mkern-1.5mu)}$. 
We reproduce it below in the case of linear polynomials over $\Q$.
In this case it was already used in \cite[Section 5]{CSan}.

\medskip

{\bf Hypothesis} ${\rm (H_1\mkern-1.5mu)}$
{\em Let $e_1,\ldots,e_r$ be pairwise different rational numbers.
Let $S$ be a finite set of primes containing
the prime factors of the denominators of $e_1,\ldots,e_r$
and the primes $p\leq r$. Suppose that we are given
$\tau_p\in\Q_p$ for $p\in S$ and a positive real constant $C$.
Then there exist $\tau\in\Q$ and primes $p_1,\ldots,p_r$ not in $S$ 
such that
\begin{enumerate}
\item $\tau$ is arbitrarily close to $\tau_p$ in
the $p$-adic topology, for $p\in S$;
\item $\tau>C$;
\item $\val_p(\tau-e_i) = 0$ for any $p\notin S\cup\{p_i\}$, $i=1,\ldots,r$;
\item $\val_{p_i}(\tau-e_i)=1$ for any $i=1,\ldots,r$.
\end{enumerate}
}
\medskip

Hypothesis~${\rm (H_1\mkern-1.5mu)}$ is usually supplemented with 
the following statement.
Let $K/\Q$ be a cyclic extension unramified outside $S$.
Assuming the conclusion of~$({\rm H}_1)$ we have the following 
implication: 

{\em if $\sum_{p\in S}\inv_p(K/\Q,\tau_p-e_i)=0$ for some $i$, then
$p_i$ splits completely in $K/\Q$.}

Hypothesis~$({\rm H}_1)$ and its supplement can be compared to the
following consequence of Proposition \ref{p1}.

\begin{prop}  \label{p2}
Let $e_1,\ldots,e_r$ be
rational numbers.
Let $S$ be a finite set of primes containing
the prime factors of the denominators of $e_1,\ldots,e_r$.
Suppose that we are given $\tau_p\in\Q_p$ for $p\in S$ and 
a positive real constant $C$. Then there exist
$\tau\in\Q$ and primes $p_1,\ldots,p_r$ not in $S$, with
$p_i=p_j$ if and only if $e_i=e_j$, such that
\begin{enumerate}
\item $\tau$ is arbitrarily close to $\tau_p$ in
the $p$-adic topology, for $p\in S$;
\item $\tau>C$;
\item $\val_p(\tau-e_i)\leq 0$ for any $p\notin S\cup\{p_i\}$, $i=1,\ldots,r$;
\item $\val_{p_i}(\tau-e_i)=1$ for any $i=1,\ldots,r$;
\item for any cyclic extension $K/\Q$ unramified outside $S$ and
such that $$\sum_{p\in S}\inv_p(K/\Q,\tau_p-e_i)=c\in\Q/\Z$$ for some $i$,
we have $\inv_{p_i}(K/\Q,\tau-e_i)=-c$; in particular, if $c=0$, then
$p_i$ splits completely in $K/\Q$.
\end{enumerate}
\end{prop}
{\em Proof.} By increasing the list of~$e_i$'s we may assume that
$r\geq 2$ and $e_i\not=e_j$ for some $i\not=j$.
Let us also assume $C>e_i$ for $i=1,\ldots,r$.
We then apply Proposition \ref{p1} to
$(\lambda_p,\mu_p)=(\tau_p,1)$ for $p\in S$.
This produces $(\lambda,\mu)\in\Z_S^2$ such that $\tau=\lambda/\mu$
satisfies all the properties in the proposition.
Indeed, (1) and (2) are clear. For $p\notin S$ we have $\val_p(\mu)\geq 0$.
For $p\notin S\cup\{p_i\}$ we have $\val_p(\lambda-e_i\mu)=0$, so that
$$\val_p(\tau-e_i)=\val_p(\lambda-e_i\mu)-\val_p(\mu)\leq 0,$$
which proves (3). We claim that $\val_{p_i}(\mu)=0$
for $i=1,\ldots,r$. Indeed, $\val_{p_i}(\mu)> 0$ implies
$\val_{p_i}(\lambda)>0$; taking~$j$ such that $e_i\not=e_j$, we obtain
$\val_{p_i}(\lambda-e_j\mu)>0$,
thus contradicting property (3) of Proposition \ref{p1}.
This proves (4). Since
$(\lambda,\mu)$ is close to $(\tau_p,1)$ in the $p$-adic topology
for $p\in S$, by continuity we have
$$\sum_{p\in S}\inv_p(K/\Q,\lambda-e_i\mu)=c.$$
We also have $\lambda-e_i\mu>0$, hence $\inv_\R(K/\Q,\lambda-e_i\mu)=0$.
By the global reciprocity law
of class field theory this implies $$\sum_{p\notin S}\inv_p(K/\Q,\lambda-e_i\mu)=-c.$$
Since $K/\Q$ is unramified outside $S$, we have $\inv_p(K/\Q,\lambda-e_i\mu)=0$
for any prime $p\notin S\cup\{p_i\}$, because in this case
$\val_p(\lambda-e_i\mu)=0$. Thus $$\inv_{p_i}(K/\Q,\lambda-e_i\mu)=-c.$$
In the case when $c=0$, we deduce from
$\val_{p_i}(\lambda-e_i\mu)=1$ that $p_i$
splits completely in $K$, so that (5) is proved. $\Box$

\medskip

\noindent{\bf Remark}. One can give a stronger variant of this
proposition by proving that $\tau$, in addition to properties (1) to (5),
can be chosen in a given Hilbertian subset of~$\Q$. The proof uses
Ekedahl's effective version of Hilbert's irreducibility theorem and the fact that
any Hilbertian subset of $\Q$ is open in the topology induced by the 
product topology of $\prod_p\Q_p$, where the product is over all primes.
We do not give a detailed proof, because this variant will not be used in this paper.

\section{Varieties fibred over the projective line} \label{fib}

\subsection{Main results} 
Recall that for a variety $X$ over $\Q$ one denotes
the image of $\Br(\Q)$ in $\Br(X)$ by $\Br_0(X)$.
If $\pi:X\to \P^1$ is a dominant
morphism of integral varieties over $\Q$, then
the corresponding vertical Brauer group is defined as follows:
$$\Brv(X)=\Br(X) \cap \pi^*\Br(\Q(\P^1))\subset \Br(\Q(X)).$$
By a $\Q$-fibre of $\pi:X\to \P^1$ we understand a fibre
above a $\Q$-point of $\P^1$.

We denote the completions of $\Q$ by $\Q_v$, and
the ring of ad\`eles of $\Q$ by $\AA_\Q$.
Let $\ov \Q$ be an algebraic closure of $\Q$. For a subfield $K\subset\ov\Q$
we write $\Gamma_K$ for the Galois group ${\rm Gal}(\ov \Q/K)$.

\begin{theorem} \label{t2}
Let $X$ be a geometrically integral variety over $\Q$
with a smooth and surjective morphism $\pi:X\to \P^1$ such that
\smallskip

{\rm (a)} with the exception of finitely many $\Q$-fibres,
denoted by $X_1,\ldots,X_r$, 
each fibre of $\pi$ contains a geometrically
integral irreducible component;

\smallskip
{\rm (b)} for each~$i$, the fibre $X_i$ contains an irreducible
component such that the algebraic closure of
$\Q$ in its field of functions is an abelian extension of $\Q$.
\smallskip

\noindent Then $\P^1(\Q)\cap \pi(X(\AA_\Q))$ is dense in
$\pi\big(X(\AA_\Q)^\Brv\big)\subset \P^1(\AA_\Q)=\prod_v\P^1(\Q_v)$.
\end{theorem}

Note that the assumptions of Theorem~\ref{t2} imply that
the generic fibre of $\pi:X\to\P^1$ is geometrically integral.
Thus all but finitely many $\Q$-fibres of $\pi$ are geometrically
integral. The cokernel of the natural map $\Br(\Q)\to\Brv(X)$
is finite by \cite[Lemma~3.1]{CTSk}.
Note finally that when $r=1$, the statement of Theorem~\ref{t2} is 
well known, see Section 2.1 of \cite{CSS}.

\medskip

{\em Proof of Theorem~\ref{t2}.}
Without loss of generality we can assume that
$X_i$ is the fibre above a point $e_i\in\A^1(\Q)$, for
$i=1,\ldots,r$. Let $K_i$ be the abelian extension of $\Q$ as in (b).

We follow the proof of \cite[Thm.~1.1]{CSS}
which uses the same method as \cite[Thm.~4.2]{CS}.
Steps~1 and~2 of the proof below repeat
the proof of \cite[Thm.~1.1]{CSS} {\em verbatim}.
Step~3 contains modifications necessary to replace the use of
Hypothesis~${\rm (H_1\mkern-1.5mu)}$ by Proposition \ref{p2}.

{\em Step} 1. Let us recall a well-known description of $\Brv(X)$.
Write $K_i$ as a compositum of cyclic extensions $K_{ij}/\Q$,
and let $\chi_{ij}:\Gamma_\Q\to \Q/\Z$ be a character such that $K_{ij}$
is isomorphic to the invariant subfield of ${\rm Ker}(\chi_{ij})$.
Let $t$ be a coordinate on $\A^1\subset\P^1$ so that $\Q(\P^1)=\Q(t)$.
The class $$A_{ij}=(K_{ij}/\Q,t-e_i)\in \Br(\Q(t))$$ of the
corresponding cyclic algebra is
ramified on $\P^1$ only at $e_i$ and $\infty$ with residues
$\chi_{ij}$ and $-\chi_{ij}$, respectively.
Let $A\in \Br(\Q(t))$ be such that $\pi^*A\in \Br(X)$.
Assumptions~(a) and~(b), together with \cite[Prop.~1.1.1]{CS},
imply that $A$ on $\P^1$ is
unramified away from $e_1,\ldots,e_r$, and that the
residue of~$A$ at~$e_i$ belongs to
$${\rm Ker}[{\rm Hom}(\Gamma_\Q,\Q/\Z)\to {\rm Hom}(\Gamma_{K_i},\Q/\Z)].$$
This group is generated by the characters $\chi_{ij}$.
Hence there exist $n_{ij}\in\Z$ such that
$A-\sum n_{ij} A_{ij}$ is unramified on $\A^1$. Since $\Br(\A^1)=\Br(\Q)$
we conclude that
$A=\sum n_{ij} A_{ij}+A_0$ for some $A_0\in\Br(\Q)$, and this implies,
by considering residues at~$\infty$, that
\begin{equation}
\sum n_{ij}\chi_{ij}=0\in {\rm Hom}(\Gamma_\Q,\Q/\Z). \label{e1}
\end{equation}
Therefore, every element of $\Brv(X)$ is of
the form $\sum n_{ij} \pi^*A_{ij}+A_0$ for some $n_{ij}$
satisfying (\ref{e1}) and some $A_0\in\Br(\Q)$.

{\em Step} 2.
We can clearly assume that $X(\AA_\Q)^\Brv\not=\emptyset$,
otherwise there is nothing to prove. Pick any
$(M_p)\in X(\AA_\Q)^\Brv$, where $M_0$ is a point in $X(\R)$. 
By a small deformation we can assume that
$M_p$ does not belong to any of the fibres $X_1,...,X_r$.

We include the real place in the finite set of places $S$
where we need to approximate. The set of real points 
$M_0 \in V(\R)$ for which $\pi(M_0)\in\P^1(\Q)$ is dense in $V(\R)$, 
and so it is enough to approximate adelic points $(M_p)$ 
such that $\pi(M_0)\in\P^1(\Q)$. By a change of variables
we then assume that $\pi(M_0) = \infty$. By another small 
deformation of the points $M_p$ for each prime $p$ 
we can further assume that $\pi(M_p) \neq \infty$ when $p\not=0$.

We include in $S$ the primes of bad reduction for $X$.
We ensure that $e_i\in\Z_S$ for each $i=1,\ldots,r$,
$e_i-e_j\in\Z_S^*$ for all $i\not=j$,
and no prime outside of $S$ is ramified in any of the fields $K_i$.
Furthermore, we increase $S$ so that
if~$K_i$ has a place of degree~$1$ over $p\notin S$,
then the corresponding $\F_p$-component of the degenerate fibre
of~$\pi$ over the reduction of~$e_i$ has an $\F_p$\nobreakdash-point.
This is achieved by using the Lang--Weil estimate, see \cite[Lemma~1.2]{CSS}.
By a similar argument we assume that on the reduction of
$X$ modulo $p\notin S$ any geometrically integral component
of a fibre over an $\F_p$-point contains an $\F_p$-point. 
All these $\F_p$\nobreakdash-points are
smooth, because $\pi$ is a smooth morphism.

Since $X(\AA_\Q)^\Brv\not=\emptyset$, by the result of Step~1 we can use
Harari's `formal lemma' \cite[Cor.~2.6.1]{H} to increase $S\subset S_1$
and choose $M_p\in X(\Q_p)$ for $p\in S_1\setminus S$ 
away from the fibres $X_1,\ldots,X_r$ so that for all $i,j$ we have
\begin{equation}
\sum_{p\in S_1} \inv_p \big(A_{ij}(\pi(M_p))\big)=0. \label{e2}
\end{equation}

{\em Step} 3. Let $\tau_p$ be the coordinate of $\pi(M_p)$,
where $p$ is a prime in $S_1$. An application of
Proposition \ref{p2} produces $\tau\in\Q$ which is
an arbitrarily large positive real number,
and is close to $\tau_p$ in the $p$-adic topology for the primes $p\in S_1$.

Let us prove that $X_\tau(\AA_\Q)\not=\emptyset$.
By the inverse function theorem we have $X_\tau(\R)\not=\emptyset$ and
$X_\tau(\Q_p)\not=\emptyset$ for $p\in S_1$. Thus it remains to consider the following two cases.

$\Q_v=\Q_p$, where $p=p_i$, $i=1,\ldots,r$.
Since $\val_{p_i}(\tau-e_i)=1$, the reduction of $\tau$ modulo~$p_i$
equals the reduction of $e_i$. In view of (\ref{e2}) property (5) of Proposition \ref{p2}
implies that for each given value of $i$ all the cyclic fields $K_{ij}$
are split at $p_i$, and thus $K_i$ is also split. Hence
there is a geometrically integral irreducible component of the $\F_{p_i}$-fibre
over the reduction of $e_i$ modulo~$p_i$. We arranged that it has
an $\F_{p_i}$-point. By Hensel's lemma it gives rise to
a $\Q_{p_i}$-point in $X_\tau$.

$\Q_v=\Q_p$, where $p\notin S_1\cup\{p_1,\ldots,p_r\}$.
We have $\val_{p}(\tau-e_i)\leq 0$ for each $i=1,\ldots,r$, and hence
the reduction of $\tau$ modulo $p$
is a point of $\P^1(\F_p)$ other than the reduction of any of
$e_1,\ldots,e_r$. Thus any $\F_p$-point
on a geometrically integral irreducible component of the fibre
at $\tau \bmod p$ gives rise
to a $\Q_p$-point on $X_\tau$, by
Hensel's lemma.

In both cases we constructed a $\Q_p$-point that comes from a $\Z_{p}$-point
on an integral model of $X_\tau$, therefore $X_\tau(\AA_\Q)\not=
\emptyset$. The theorem is proved. $\Box$

\begin{cor} \label{c1}
In the situation of Theorem~\ref{t2}, let us further assume that
all but finitely many
$\Q$-fibres of $\pi:X\to\P^1$ satisfy the Hasse principle. Then
$\pi(X(\Q))$ is dense in $\pi\big(X(\AA_\Q)^\Brv\big)$.
If, in addition, these $\Q$-fibres $X_\tau$ are such that
$X_\tau(\Q)$ is dense in $X_\tau(\AA_\Q)$, then $X(\Q)$ is dense in $X(\AA_\Q)^\Brv$.
\end{cor}

\noindent{\bf Remark}.
If the generic fibre of $\pi:X\to\P^1$ is proper, then all
but finitely many fibres of $\pi$ are proper. For proper $\Q$-fibres $X_\tau$
the approximation assumption in Corollary \ref{c1}
is that of {\em weak approximation}, since in this case
$X_\tau(\AA_\Q)=\prod_v X_\tau(\Q_v)$.
By Hironaka's theorem one can always replace \mbox{$\pi:X\to\P^1$} by a partial
compactification $\pi':X'\to\P^1$ such that $X$ is a dense open subset
of $X'$ and the morphism $\pi'$ is smooth with proper generic fibre.

\medskip

We now give a statement for a smooth and proper variety $X$,
to be compared with \cite[Thm.~1.1]{CSS}.

\begin{theorem} \label{t3}
Let $X$ be a smooth, proper and geometrically integral variety
over $\Q$ with a surjective morphism $\pi:X\to\P^1$ such that
\smallskip

{\rm (a)} with the exception of finitely many $\Q$-fibres,
denoted by $X_1,\ldots,X_r$, each fibre of $\pi$ contains a geometrically
integral irreducible component of multiplicity one;

\smallskip
{\rm (b)} for each~$i$, the fibre $X_i$ contains an irreducible
component of multiplicity one such that the algebraic closure of
$\Q$ in its field of functions is an abelian extension of $\Q$.
\smallskip

\noindent Then $\P^1(\Q)\cap \pi(X(\AA_\Q))$ is dense in
$\pi\big(X(\AA_\Q)^\Brv\big)\subset \P^1(\AA_\Q)=\prod_v\P^1(\Q_v)$.
If, moreover, all but finitely many $\Q$-fibres of $\pi$ satisfy
the Hasse principle and weak approximation, then $X(\Q)$ is dense in $X(\AA_\Q)^\Brv$.
\end{theorem}

{\em Proof.} Let $Y \subset X$ denote the smooth locus of~$\pi$.
By (a) and (b) each fibre of $\pi:X\to\P^1$ contains
a multiplicity one irreducible component, hence $\pi(Y)=\P^1$.
Thus Theorem~\ref{t2} can be applied to $\pi:Y\to\P^1$.
We have $\Brv(X)=\Br(X)\cap\Brv(Y)$, and by \cite[Lemma~3.1]{CTSk}, the set
$Y(\AA_\Q)^{\Brv(Y)}$ is dense in $X(\AA_\Q)^{\Brv(X)}$.
This proves the first statement. The second one now follows from
Corollary \ref{c1}. $\Box$

\medskip

\noindent{\bf Remarks}. (1) We note that assumptions (a) and (b)
will hold for any smooth, proper and geometrically integral variety
$X'$ over $\Q$ with a surjective morphism $\pi':X'\to\P^1$
such that the generic fibres of $\pi$ and $\pi'$ are isomorphic.
This follows from \cite[Lemma~1.1]{S96}, see also \cite[Lemme 3.8]{W}.

(2) Using the remark after Proposition \ref{p2} one can prove
a stronger variant of Theorem \ref{t3} with ``all but finitely many
$\Q$-fibres of $\pi$" replaced by ``the $\Q$-fibres of $\pi$ over
a Hilbertian subset of $\Q$". The details are left to the interested reader.

\subsection{Application to pencils of Severi--Brauer and similar varieties}
Theorem~\ref{t3} can be applied to the
varieties considered by Colliot-Th\'el\`ene
and Swinnerton-Dyer in \cite{CS}, see also \cite{CSS}.

\begin{cor} \label{cc1}
Let $X$ be a smooth, proper and geometrically integral variety over $\Q$
with a morphism $\pi:X\to\P^1$. Suppose that the generic fibre of $\pi$
is a Severi--Brauer variety (for example, a conic),
a $2$-dimensional quadric, or a product of such.
If all the fibres of $\pi$ that are not geometrically integral are
$\Q$-fibres, then $X(\Q)$ is dense in $X(\AA_\Q)^\Brv$.
\end{cor}

{\em Proof.} The assumptions of Theorem~\ref{t3} are satisfied
for $X$ and $\pi$ by well-known results of class field theory
and by the structure of fibres of regular models of quadric
bundles \cite{S90} and Artin models of Severi--Brauer
varieties \cite{F} (or see \cite[pp.~117--118]{W} for
an alternative argument).  $\Box$

\medskip

In particular, this result gives a uniform approach to
Theorems 1.1, 1.3 and~1.4 of \cite{BMS} that does not use
descent and leads to shorter and more natural proofs.
See \cite{BMS} for a survey
of known results on conic and quadric bundles, most of which are established
over an arbitrary number field in place of $\Q$ but under
a strong restriction on the number of degenerate geometric fibres.

\subsection{Theorem~\ref{t1} as a consequence of Theorem~\ref{t2}} \label{nn}
Here is one way to deduce Theorem~\ref{t1} from Corollary~\ref{c1}
that keeps the Brauer group calculations to the minimum.

Let $W$ be the quasi-affine subvariety of
$\A^2\times \A^{d_1}\times\ldots\times\A^{d_r}$ given by
$$b_i(u-e_iv)=N_{K_i/\Q}(\bfx_i), \quad i=1,\ldots,r, \quad (u,v)\not=(0,0).$$
The variety $V$ defined in (\ref{E1}) is a dense open subset of $W$.
The projection to the coordinates $(u,v)$ defines a morphism
$W\to\A^2\setminus (0,0)$. Let $\pi: W\to\P^1$ be the composed morphism
$W\to\A^2\setminus (0,0)\to\P^1$, and let $X\subset W$
be the smooth locus of $\pi$.
It is easy to see that $\pi(X)=\P^1$. Let $\pi':Y\to \P^1$
be a partial compactification of $\pi:X\to \P^1$ as in the remark 
after Corollary~\ref{c1}, so that $\pi'$ is smooth with proper 
generic fibre.

Let $t=u/v$ be a coordinate on $\P^1$. It is straightforward to see
that conditions (a) and (b) of Theorem~\ref{t2} hold. In order to deduce
Theorem~\ref{t1} from Corollary~\ref{c1} we need to prove that
\begin{enumerate}
\item geometrically integral, proper $\Q$-fibres of $\pi'$
satisfy the Hasse principle and weak approximation;
\item $\Brv(Y)=\Br_0(Y)$.
\end{enumerate}
The fibre of $\pi$ at $\tau\in\Q$, $\tau\not=e_i$,
is the affine variety
$$\frac{N_{K_1/\Q}(\bfx_1)}{b_1(\tau-e_1)}=\ldots=\frac{N_{K_r/\Q}(\bfx_r)}{b_r(\tau-e_r)}\not=0.$$
This is a principal homogeneous space of the torus $T$ defined by
$$N_{K_1/\Q}(\bft_1)=\ldots=N_{K_r/\Q}(\bft_r)\not=0.$$
The Hasse principle and weak approximation hold
for smooth compactifications of principal homogeneous spaces of $T$
if ${\cyr X}^2_\omega(\Q,\hat T)=0$, see \cite[Ch.~8]{Sansuc}
or \cite[Thm.~6.3.1]{Sk}.
Let $T_i=R^1_{K_i/\Q}(\G_m)$ be the norm torus attached to 
$K_i/\Q$, and let $G_i$ be the cyclic group ${\rm Gal}(K_i/\Q)$.
Then $T$ is an extension of $\G_m$ by the product of the tori $T_i$, so that
we have the exact sequence of $\Gamma_\Q$-modules of characters
\begin{equation}
0\to\Z\to\hat T\to \prod_{i=1}^r \hat T_i \to 0.
\label{E2}
\end{equation}
The long exact sequence of Galois cohomology
gives rise to the exact sequence
$$\prod_{i=1}^r {\rm Hom}(G_i,\Q/\Z) \to {\rm Hom}(\Gamma_\Q,\Q/\Z)\to
H^2(\Q,\hat T)\to
\prod_{i=1}^r H^2(\Q,\hat T_i).$$
Let $K$ be the compositum of the fields $K_i$, and let
$G={\rm Gal}(K/\Q)$. 
Since ${\cyr X}^2_\omega(\Q,\hat T_i)=0$ it is enough to show that
if $\alpha\in {\rm Hom}(\Gamma_\Q,\Q/\Z)$ goes to
${\cyr X}^2_\omega(\Q,\hat T)$ in $H^2(\Q,\hat T)$, then $\alpha
\in {\rm Hom}(G,\Q/\Z)$. But the tori $T_i$ and $T$ are split by $K$,
hence (\ref{E2}) is split as an extension of
$\Gamma_K$-modules. It follows that the restriction of $\alpha$
to $H^2(K,\Z)={\rm Hom}(\Gamma_K,\Q/\Z)$ is in
${\cyr X}^2_\omega(K,\Z)=0$, thus $\alpha\in {\rm Hom}(G,\Q/\Z)$.
This proves (1).

Since $\pi:X\to \P^1$ factors through the inclusion of $X$ into $Y$,
to prove (2) it is enough to prove $\Brv(X)=\Br_0(X)$.
Write $\chi_i\in {\rm Hom}(\Gamma_\Q,\Q/\Z)$ for the image of
a generator of ${\rm Hom}(G_i,\Q/\Z)$.
Recall from Step~1 of the proof of Theorem~\ref{t2}
that any $A\in\Br(\Q(t))$ such that $\pi^*A\in\Br(X)$
can be written as
$$A=\sum_{i=1}^r n_i (\chi_i,t-e_i)+A_0, \quad \text{where}
\quad \sum_{i=1}^r n_i\chi_i=0,$$
for some $A_0\in\Br(\Q)$. In $\Br(X)$ we have
$$\pi^*A=\sum_{i=1}^r n_i (\chi_i,u-e_iv)-\sum_{i=1}^r n_i (\chi_i,v)+A_0=
-\sum_{i=1}^r n_i (\chi_i,b_i)+A_0\in\Br_0(X),$$
since $(\chi_i,N_{K_i/\Q}(\bfx_i))=0$ in $\Br(\Q(X))$.
This finishes the proof of (1) and (2).

\section{Norms and their products as products of linear polynomials} \label{four}

\subsection{Cyclic extensions}
Consider the following system of Diophantine equations:
\begin{equation}
N_{K_i/\Q}(\bfx_i)=P_i(t),\quad i=1,\ldots,r,\label{eq1}
\end{equation}
where $K_i/\Q$ are cyclic extensions and the polynomials
$P_i(t)$ are products of (possibly repeated) linear factors over $\Q$.

\begin{cor} \label{cc2}
Let $X$ be a smooth, proper and geometrically integral variety over $\Q$
with a morphism $\pi:X\to \P^1$ such that the generic fibre
of $\pi$ is birationally equivalent
to the affine variety {\rm (\ref{eq1})} over $\Q(\P^1)=\Q(t)$.
Then $X(\Q)$ is dense in $X(\AA_\Q)^\Brv$.
\end{cor}

{\em Proof.} Each fibre of $\pi$
outside infinity and the zero set of $P_1(t)\ldots P_r(t)=0$
contains a geometrically integral irreducible component of multiplicity
one. Since $\pi$ has a section over the compositum $K_1\ldots K_r$, which
is an abelian extension of $\Q$,
assumptions (a) and (b) of Theorem~\ref{t3} hold for $\pi:X\to\P^1$.
By Hasse's norm theorem the varieties $N_{K/\Q}(\bfz)=c$,
where $K/\Q$ is cyclic and $c\in\Q^*$, satisfy the Hasse principle.
Moreover,
smooth and proper models of principal homogeneous spaces
of cyclic norm tori satisfy the Hasse principle and weak approximation,
by \cite[Ch.~8]{Sansuc}. We conclude by Theorem~\ref{t3}. $\Box$

\medskip

For $r=1$ the statement of Corollary~\ref{cc2} was previously
known for any finite extension $K/\Q$ in the case when $P_1$ 
has at most two roots, see \cite{HS}, \cite{CHS}; see also
\cite{SJ}, where $\Q$ was replaced by an arbitrary number field.

For any cyclic extension of fields $K/k$ the affine variety $N_{K/k}(\bfx)=c$,
where
$c\in k^*$, is well known to be birationally equivalent to the Severi--Brauer
variety defined by the cyclic algebra $(K/k,c)$. Thus
Corollary \ref{cc2} can be seen as a particular case of Corollary \ref{cc1}.

When each polynomial $P_i(t)$ is linear we have the following consequence
of Corollary 4.1.

\begin{cor} \label{cc3}
Let $K_i$ be a cyclic extension of $\Q$ of degree $d_i$, $i=1,\ldots,r$.
Let $b_i\in\Q^*$ and $e_i\in\Q$, $i=1,\ldots,r$.
Then the variety over $\Q$ defined by 
$$b_i(t-e_i)=N_{K_i/\Q}(\bfx_i)\not=0, \quad i=1,\ldots,r,$$
satisfies the Hasse principle and weak approximation.
\end{cor}
{\em Proof}. An easy calculation shows that this variety $X$ is smooth.
By Corollary \ref{cc2} is enough to prove
that $\Brv(X)=\Br_0(X)$.
In Step~1 of the proof of Theorem~\ref{t2} we saw
that for any $A\in\Br(\Q(t))$ such that $\pi^*A\in\Br(X)\subset\Br(\Q(X))$
there exists $A_0\in\Br(\Q)$ for which we can write
$$A=\sum_{i=1}^r n_i (K_i/\Q,t-e_i)+A_0.$$
Since $(K_i/\Q,N_{K_i/\Q}(\bfx_i))=0$ in $\Br(\Q(X))$,
the element $\pi^*A\in \Br(X)$ can be written as
$$\pi^*A=-\sum_{i=1}^r n_i (K_i/\Q,b_i)+A_0\ \in\ \Br_0(X). \qquad\Box$$

\medskip

The following statement is deduced from Corollary~\ref{cc2}
by an easy application of the fibration method in the form of \cite[Thm.~3.2.1]{H2}.
\begin{cor} \label{cccc}
Let $X$ be a smooth and proper model of the variety over $\Q$
defined by the system of equations
\begin{equation}
N_{K_i/\Q}(\bfx_i)=P_i(t_1,\ldots,t_n), \quad i=1,\ldots,r,\label{eq}
\end{equation}
where each $K_i$ is a cyclic extension of $\Q$
and each $P_i(t_1,\ldots,t_n)$ is a product of polynomials of degree~$1$ over~$\Q$.
Then $X(\Q)$ is dense in $X(\AA_\Q)^\Br$.
\end{cor}
In \cite{SS}, under a mild general position condition,
this was proved for $r=1$ and $\deg\,P_1\leq 2n$, but with
the cyclic extension $K$ of $\Q$ replaced by
any finite extension of an arbitrary number field.

\subsection{Products of norms}
Instead of a norm form of a cyclic extension of $\Q$ we can consider
a product of norm forms associated to field extensions of $\Q$
satisfying certain conditions. We start with one more application of
Theorem~\ref{t3}.

\begin{cor} \label{cc4}
Let $P(t)$ be a product of (possibly repeated) 
linear factors over~$\Q$. Let $L_1,\ldots,L_n$ be $n\geq 2$ finite field extensions
of $\Q$ such that $L_1/\Q$ is abelian and linearly disjoint from
the compositum $L_2\ldots L_n$.
Let $X$ be a smooth, proper and geometrically integral variety over $\Q$
with a morphism $\pi:X\to \P^1$ such that the generic fibre of $\pi$
is birationally equivalent to the affine variety 
\begin{equation}
N_{L_1/\Q}(\bfx_1)\ldots N_{L_n/\Q}(\bfx_n)=P(t)\label{eq5}
\end{equation}
over $\Q(\P^1)=\Q(t)$. Then $X$ satisfies the Hasse principle
and weak approximation.
\end{cor}
{\em Proof.} This proof is similar to that of Corollary \ref{cc2}.
Assumptions (a) and (b)
of Theorem~\ref{t3} are satisfied since $L_1/\Q$ is abelian.
To prove that almost all $\Q$-fibres satisfy the Hasse principle
and weak approximation it is enough, by \cite[Ch.~8]{Sansuc} 
(see also \cite[Thm.~6.3.1]{Sk}),
to verify that ${\cyr X}^2_\omega(\Q,\hat T)=0$, where $T$
is the multinorm torus over $\Q$ attached to the fields $L_1,\ldots,L_n$.
This was proved by Demarche and Wei \cite[Thm.~1]{DW}. 
To finish the proof we note that $\Brv(X)=\Br_0(X)$.
Indeed, let $A\in \Br(\Q(t))$ be such that $\pi^*A\in\Br(X)$. 
The morphism $\pi$ has a section defined over $L_i$,
for each $i=1,\ldots,n$.
By restricting $A$ to this section we see the restriction of 
$A$ to $\Br(L_i(t))$ comes from
$\Br(\P^1_{L_i})=\Br(L_i)$. In particular, the residues
of $A$ at the roots of $P(t)$ are in the kernel of the
map $H^1(\Q,\Q/\Z)\lra H^1(L_i,\Q/\Z)$. Since  
$L_1\cap L_2\ldots L_n=\Q$ there is no non-trivial cyclic extension of $\Q$
contained in all of the $L_i$. This implies that $A$
is not ramified at the zero set of $P(t)$. A well known
argument \cite[Prop.~1.1.1]{CS} shows that $A$ is unramified away
from the zero set of $P(t)$. Hence $A\in\Br(\A^1)=\Br(\Q)$. $\Box$
\medskip

For more cases when ${\cyr X}^2_\omega(\Q,\hat T)=0$ for the multinorm
torus $T$ see \cite{DW}, Thm.~1 and Cor.~8.
One can extend Corollary \ref{cc4} to systems of equations (\ref{eq5})
and replace $P(t)$ by a product of linear
polynomials in several variables. We leave the details to
the interested reader.

\medskip

Following Wei \cite[Thm.~3.5]{Wei} we now
consider a case where the $\Q$-fibres do not satisfy
the Hasse principle nor weak approximation. 

\begin{prop} \label{pp1}
Let $P(t)$ be a product of (possibly repeated) linear factors over $\Q$, and let
$a,b\in\Q^*$. 
Let $X$ be a smooth, proper and geometrically integral variety over $\Q$
with a morphism $\pi:X\to \P^1$ such that the generic fibre of $\pi$
is birationally equivalent to the affine variety 
\begin{equation}
N_{\Q(\sqrt{a})/\Q}(\bfx) N_{\Q(\sqrt{b})/\Q}(\bfy) N_{\Q(\sqrt{ab})/\Q}(\bfz)=P(t)
\label{eq6}
\end{equation}
over $\Q(\P^1)=\Q(t)$. Then $X(\Q)$ is dense in $X(\AA_\Q)^\Br$.
\end{prop}
{\em Proof.} We can assume that $\Q(\sqrt{a})$, $\Q(\sqrt{b})$ and $\Q(\sqrt{ab})$
are quadratic fields, otherwise the variety $X$ is rational and the statement
is clear. Let $V$ be the smooth locus of the affine variety 
(\ref{eq6}), and let $U$ be the image of $V$ by the projection 
to the coordinate $t$. It is clear that
$\P^1\setminus U$ is a finite union of $\Q$-points. 
The fibres of $V\to U$ are principal homogeneous spaces of 
the torus $T$ that is given by
$$N_{\Q(\sqrt{a})/\Q}(\bfx)N_{\Q(\sqrt{b})/\Q}(\bfy)N_{\Q(\sqrt{ab})/\Q}(\bfz)=1.$$
Let $E$ be a smooth equivariant compactification of $T$
(which exists by \cite{CHS1}), and let $V^c=V\times^T E$
be the contracted product. Then $V^c\to U$ is a fibre-wise smooth 
compactification of 
$V\to U$. We take $\pi:X\to\P^1$ such that $X\times_{\P^1}U=V^c$.
We compose $\pi$ with an automorphism of $\P^1$ to ensure that the fibre at infinity 
is smooth and is close to the real point that we need to approximate;
in particular, the fibre at infinity contains a real point.
An obvious change of variables shows that 
$X$ contains an open set which is the smooth locus of the affine variety
given by 
$$N_{\Q(\sqrt{a})/\Q}(\bfx)N_{\Q(\sqrt{b})/\Q}(\bfy)N_{\Q(\sqrt{ab})/\Q}(\bfz)=Q(t),
$$ where $Q(t)$ is a polynomial with rational roots $e_1,\ldots, e_r$ 
such that $U$ is the complement to $\{e_1,\ldots, e_r\}$ in $\P^1$.
Note that for any $\tau\in U(\Q)$ we have
$X_\tau(\AA_\Q)\not=\emptyset$ by \cite[Prop.~5.1]{C}.

The quaternion algebra 
$A=(N_{\Q(\sqrt{a})/\Q}(\bfx),b)$ defines an element of $\Br(\pi^{-1}(U))$. 

We are given points $M_p\in X(\Q_p)$ for all primes $p$
and $M_0\in X(\R)$ such that $(M_p)\in X(\AA_\Q)^\Br$. 
Since $\Br(X)/\Br_0(X)$ is finite, by
a small deformation we can assume that $\pi(M_p)$ is a point
in $U\cap\A^1$ where $t$ equals $\tau_p\in \Q_p$.

Let $S_0$ be the finite set of places of $\Q$ where we need to approximate.
We can find a finite set $S$ of places containing $S_0$ and the real place,
such that $\pi: X\to\P^1$ extends to a proper
morphism $\pi:\calX\to \P^1_{\Z_S}$ with $\calX$ regular.
By doing so we can ensure that
$S$ contains the primes where at least one 
of our quadratic fields is ramified, and that we have $a,\,b\in\Z_S$,
$Q(t)\in\Z_S[t]$, and $e_i\in\Z_S$ for $i=1,\ldots, r$.
By Harari's `formal lemma' \cite[Cor.~2.6.1]{H} we can further enlarge $S$ so that 
$$\sum_{p\in S} \inv_p(A(M_p))=0, \quad 
\sum_{p\in S} \inv_p(b,\tau_p-e_i)=0, \ i=1,\ldots, r.$$
(For this we may need to modify the points $M_p$ for $p\in S\setminus S_0$.)
Let $\calU$ be the complement to the Zariski closure of 
$e_1\cup\ldots\cup e_r$ in $\P^1_{\Z_S}$. The 
same algebra $A$ defines a class in $\Br(\pi^{-1}(\calU))$.
An application of Proposition \ref{p2} gives a $\Q$-point $\tau$
in $U\cap \A^1$ that is as large as we want in the real topology and is
close to~$\tau_p$ in the $p$-adic topology for the primes $p\in S$. 
For $p\notin S\cup\{p_1,\ldots,p_r\}$ we see from property (3) of 
Proposition \ref{p2} that the Zariski closure of $\tau$ in 
$\P^1_{\Z_p}$ is contained in $\calU\times_{\Z_S}\Z_p$. 
This implies that for any $N_p\in X_\tau(\Q_p)$
the value $A(N_p)\in\Br(\Q_p)$ comes from $\Br(\Z_p)=0$.
From property (5) we see that for
each $i=1,\ldots,r$ the prime $p_i$ splits in $\Q(\sqrt{b})$, hence 
$A(N_{p_i})=0$ for any $N_{p_i}\in X_\tau(\Q_{p_i})$.
By continuity and the inverse function theorem
we can find $N_p\in X_\tau(\Q_p)$ close enough 
to $M_p$, for $p\in S$ (including the real place), so that 
$\sum_{p\in S} \inv_p(A(N_p))=0$. Summing over all places of $\Q$
we now have $\sum_{p} \inv_p(A(N_p))=0$, for any choice
of $N_p$, $p\notin S$. By \cite[Thm.~4.1]{C} the algebra
$A$ generates $\Br(X_\tau)$ modulo the image of $\Br(\Q)$.
By \cite[Ch.~8]{Sansuc} or \cite[Thm.~6.3.1]{Sk} the set
$X_\tau(\Q)$ is dense in $X_\tau(\AA_\Q)^\Br$, so 
we can find a $\Q$-point in $X_\tau$ close to $M_p$, for
$p\in S$. $\Box$

\subsection{Non-cyclic extensions of prime degree}
The method of Colliot-Th\'el\`ene and Wei
\cite[Thm.~3.6]{Wei} can be used to prove the following result.

\begin{theorem} \label{pp2}
Let $P(t)$ be a product of (possibly repeated) 
linear factors over $\Q$. Let $K$ be a non-cyclic extension
of $\Q$ of prime degree such that the Galois group
of the normal closure of $K$ over $\Q$ has a non-trivial
abelian quotient.
Let $X$ be a smooth, proper and geometrically integral variety over $\Q$
with a morphism $\pi:X\to \P^1$ such that the generic fibre of $\pi$
is birationally equivalent to the affine variety 
$N_{K/\Q}(\bfx)=P(t)$ over $\Q(t)$. 
Then $X$ satisfies the Hasse principle and weak approximation.
\end{theorem}

This result covers the `generic' case of the field $K=\Q[t]/(f(t))$, 
for which $f(t)$ is a polynomial of prime degree $\ell>2$ 
such that the Galois group
of $f(t)$ is the symmetric group $S_\ell$ (or the dihedral group $D_{2\ell}$.)
In particular, the assumption on $K$ in Theorem \ref{pp2} is 
automatically satisfied when $[K:\Q]=3$.

\medskip

{\em Proof of Theorem \ref{pp2}.} 
We start in the same way as in the proof of Proposition \ref{pp1}.
We can assume that
$X$ contains an open set which is the smooth locus of the affine variety
$N_{K/\Q}(\bfx)=Q(t)$, where $Q(t)$ is a product of powers of
$t-e_i$, $i=1,\ldots, e_r$,
with the additional assumption that the fibre at infinity is smooth and
contains a real point close to the real point that we want to approximate.

\begin{lemma} \label{le3}
We have $\Br(X)=\Br_0(X)$.
\end{lemma}
{\em Proof.} Let $T$ be the norm torus $N_{K/\Q}(\bfx)=1$.
Since $\ell=[K:\Q]$ is prime, \cite[Prop.~9.1, Prop.~9.5]{CSan1} 
gives $H^1(F,\Pic(Z\times_F\ov F))={\cyr X}^2_\omega(F,\hat T)=0$
for any smooth and proper variety $Z$ over a field $F$ such that
a dense open subset of $Z$ is a principal homogeneous space of $T$.
Applying this to the generic fibre of $\pi:X\to\P^1$ we see
that $\Br(X)=\Brv(X)$. 

Now let $A\in \Br(\Q(t))$ be such that $\pi^*A\in\Br(X)$. 
The morphism $\pi$ has a section defined over $K$. 
By restricting to it we see that the image of 
$A$ in $\Br(K(t))$ belongs to the injective image of
$\Br(\P^1_K)=\Br(K)$. In particular, the residue
of $A$ at $e_i$ lies in the kernel of the
map $$H^1(\Q,\Q/\Z)\lra H^1(K,\Q/\Z).$$ Since $K$ contains no
cyclic extension of $\Q$, this kernel is zero. Thus $A$
is not ramified at the zero set of $Q(t)$. Since $A$ is also unramified outside of
the zero set of $Q(t)$, we see that $A\in\Br(\Q)$. $\Box$

\medskip

Let $L$ be the normal closure of $K/\Q$.
By assumption there exists a cyclic extension $k/\Q$ of prime
degree such that $k\subset L$. Let $\ell=[K:\Q]$, $q=[k:\Q]$.
Since $\Gal(L/\Q)\subset S_\ell$ and $k\not=K$ we see that $q<\ell$. 

\begin{lemma} \label{le2}
Let $a\in \Q^*$. If~$p$ is a prime unramified in~$L$ and inert in~$k$,
then the equation $N_{K/\Q}(\bfx)=a$ is solvable in $\Q_p$.
\end{lemma}
{\em Proof.} Write $K\otimes_{\Q}\Q_p=K_{v_1}\oplus\ldots\oplus K_{v_s}$,
and let $d_i=[K_{v_i}:\Q_p]$. 

If $s>1$, then since $\ell=d_1+...+d_s$ 
is a prime number, there exist integers $n_1,\ldots,n_s$ such that
$1=n_1d_1+...+n_sd_s$. It follows that
$$a=\prod_{i=1}^s N_{K_{v_i}/\Q_p}(a^{n_i}) \ \in \ N_{K/\Q}(K\otimes_{\Q}\Q_p),$$
so we are done.

If $s=1$, then $K\otimes_{\Q}\Q_p=K_v$ is a field extension of $\Q_p$
of degree $\ell$. By assumption~$p$ is inert in~$k$, so that
$k\otimes_{\Q}\Q_p=k_w$ is a field. Since $[k_w:\Q_p]=q$ is a 
prime less than $\ell$, the fields $k_w$ and $K_w$
are linearly disjoint over $\Q_p$, so that~$k_wK_v$ is a field. Thus
$p$ is inert in $kK\subset L$, which implies that the
Frobenius at $p$ in $\Gal(L/\Q)$ is an element of order divisible
by $\ell q$. However, $S_\ell$ contains no such elements,
so the case $s=1$ is impossible. $\Box$

\medskip

{\em End of proof of Theorem \ref{pp2}.} 
We are given points $M_p\in X(\Q_p)$ for all primes $p$
and $M_0\in X(\R)$. By
a small deformation we can assume that $\pi(M_p)$ is a point
in $U\cap\A^1$ where $t$ equals $\tau_p\in \Q_p$.
Let $S$ be the finite set of places of $\Q$ where we need to approximate,
containing the real place and the primes of bad 
reduction for $X$. We also assume that $L$ is unramified away from $S$.
Consider the cyclic algebras 
$$A_i=(k/\Q,t-e_i)\in \Br(\Q(X)), \ i=1,\ldots,r.$$
Harari's `formal lemma' \cite[Cor.~2.6.1]{H} and Lemma \ref{le3} imply
that we can introduce new primes into $S$ and choose the 
corresponding points $M_p$ so that 
$$\sum_{p\in S}\inv_p(A_i(\tau_p))\not=0, \quad i=1,\ldots,r.$$
An application of Proposition \ref{p2} gives a $\Q$-point $\tau$
in $U\cap \A^1$ that is as large as we want in the real topology and is
close to $\tau_p$ in the $p$-adic topology for the primes $p\in S$. 
This ensures that $X_\tau(\R)\not=\emptyset$ and $X_\tau(\Q_p)\not=\emptyset$
for all $p\in S$.
For $p\notin S\cup\{p_1,\ldots,p_r\}$ we see from property (3) of 
Proposition \ref{p2} that $\tau$ reduces modulo $p$ to a point
of $\P^1(\F_p)$ other than the reduction of any of $e_1,\ldots,e_r$.
The corresponding fibre over $\F_p$ contains a principal homogeneous
space of a torus over a finite field, and hence an $\F_p$-point,
by Lang's theorem. By Hensel's lemma it gives rise to a $\Q_p$-point
in $X_\tau$. Finally, property (5) of Proposition \ref{p2}
gives that $\inv_{p_i}(A_i(\tau))\not=0$. By
property (4) this implies that $p_i$ is inert in $k$. Now 
an application of Lemma \ref{le2} shows that $X_\tau(\Q_{p_i})
\not=\emptyset$. This holds for all $i=1,\ldots,r$ so we conclude that
$X_\tau(\AA_\Q)\not=\emptyset$. To finish the proof we note that
${\cyr X}^2_\omega(\Q,\hat T)=0$ implies
that the principal homogeneous spaces of $T$ over $\Q$ satisfy the Hasse
principle and weak approximation \cite[Ch.~8]{Sansuc}. $\Box$

\bibliographystyle{amsbracket}
\providecommand{\bysame}{\leavevmode\hbox to3em{\hrulefill}\thinspace}

\end{document}